\documentclass[10pt, conference, letterpaper]{IEEEtran}

\usepackage{multicol}
\usepackage{multirow}
\usepackage{graphicx}

\usepackage{cite}
\usepackage{url}

\usepackage{amssymb}
\usepackage{amsthm, amsmath}
\usepackage{amsbsy}
\usepackage{color}

\usepackage{subfigure}

\usepackage{cite,url}

\usepackage{dsfont}  
\usepackage{algorithmic}
\usepackage[ruled]{algorithm2e}
\SetKwInput{kwInitStep}{Initialization Step}
\SetKwInput{kwGibbsStep}{Gibbs Sampling Step}
\SetKwInput{kwTrainStep}{Training}
\SetKwInput{KwPredictStep}{Prediction}

\def\squareforqed{\IEEEQED}
\def\qed{\ifmmode\squareforqed\else{\unskip\nobreak\hfil
\penalty50\hskip1em\null\nobreak\hfil\squareforqed
\parfillskip=0pt\finalhyphendemerits=0\endgraf}\fi}

\newcommand{\ud}{\mathrm{d}}

\begin{document}

\title{Recover Fine-Grained Spatial Data \\from Coarse Aggregation}

\author{\IEEEauthorblockN{Bang Liu}\thanks{ This work was partially supported by the NSERC-CRD and NSERC-RGPIN grants.}
\IEEEauthorblockA{
University of Alberta\\
bang3@ualberta.ca}
\and
\IEEEauthorblockN{Borislav Mavrin}
\IEEEauthorblockA{
University of Alberta\\
mavrin@ualberta.ca}
\and
\IEEEauthorblockN{Linglong Kong}
\IEEEauthorblockA{
University of Alberta\\
lkong@ualberta.ca}
\and
\IEEEauthorblockN{Di Niu}
\IEEEauthorblockA{
University of Alberta\\
dniu@ualberta.ca}}

\maketitle
\thispagestyle{empty}
\pagestyle{empty}

\begin{abstract}
	In this paper, we study a new type of \emph{spatial sparse recovery} problem, that is to infer the fine-grained spatial distribution of certain density data in a region only based on the aggregate observations recorded for each of its subregions.
	One typical example of this spatial sparse recovery problem is to infer spatial distribution of cellphone activities based on aggregate mobile traffic volumes observed at sparsely scattered base stations.  
	We propose a novel \emph{Constrained Spatial Smoothing} (CSS) approach, which exploits the local continuity that exists in many types of spatial data to perform sparse recovery via finite-element methods, while enforcing the aggregated observation constraints through an innovative use of the ADMM algorithm. We also improve the approach to further utilize additional geographical attributes.
	Extensive evaluations based on a large dataset of phone call records and a demographical dataset from the city of Milan show that our approach significantly outperforms various state-of-the-art approaches, including Spatial Spline Regression (SSR).

\end{abstract}

\section{Introduction}
\label{sec:intro}

In this paper, we consider a new type of \emph{spatial sparse recovery} problem, that is to infer the fine-grained spatial distribution of certain density data in a region based on the aggregate observations recorded for each of its subregions. Such a spatial sparse recovery problem is of particular interests to many big data applications, where acquiring fine-grained spatial data involves either additional technical overhead or privacy issues. For example, a telecommunication service provider may only be able to monitor the aggregate mobile traffic on each cell tower (base station), but want to infer the fine-grained spatial distribution of cell phone activities for capacity planning, device installation, urban ecology~\cite{bcici_mobihoc15}, and the study of population density estimation~\cite{douglass2015high}. Another example is to infer the spatial distribution of population densities (e.g., voter population for a certain party) based on the aggregate population observed at sparsely scattered polling stations. A further example is for an Internet media provider or retailer, e.g., Google, Amazon, etc., to reconstruct a fine-grained geographical distribution of its users only from aggregated user counts observed at certain datacenters or points of presence (PoPs). The ability of spatial recovery from coarse aggregation will be critical in all these cases, as tracking the position of each individual may incur either technical overhead or privacy concerns.

However, the above-mentioned problem is very challenging. 
A straightforward solution is to have a patched piece-wise constant estimation by assuming the density is uniformly distributed within each subregion in which there is an aggregate observation. However, this approach gives a non-smooth piece-wise constant spatial field, which offers little value in terms of discovering hot spots. 
Moreover, the densities estimated in this way will jumps between neighboring subregions, ignoring the local continuity of the spatial densities across subregions. Such local continuity, however, is present in many spatial data, e.g., cell phone/Internet activities, which are highly dependent on underlying factors such as urban geography, area functionality, and population density, all of which are usually smoothly changing. 
Alternatively, we may use spatial smoothing techniques \cite{Sanga13} developed in statistics in the past decade to smoothen the patched estimation. 
However, nearly all existing spatial smoothing techniques \cite{Wood08, Sanga13, guillas2010bivariate} aim to reconstruct a spatial field of densities based on \emph{sampled} observations, e.g., recovering a spatial field of temperatures according to temperature readings at a few sample points, instead of based on coarse-grained \emph{aggregate observations}. As a result, applying existing spatial smoothing techniques to our new problem may violate the necessary constraint that the estimated spatial data in each subregion must sum up to its corresponding aggregate observation in the first place, leading to systematic errors.

In this paper, we propose a new approach called \emph{Constrained Spatial Smoothing} (CSS) to overcome the difficulties mentioned above. Specifically, we aim to reconstruct a spatial field $f$ of densities over a region $\Omega$ by penalizing the ``roughness'' of $f$ (especially across subregions), subject to the constraint that the aggregation of discretized values of $f$ in each patched subregion $\Omega_{B_i}$ equals to the aggregate value $z_i$ observed in $\Omega_{B_i}$. Our new approach is different from traditional spatial smoothing due to the presence of the additional constraint.
We propose an Alternating Direction Method of Multipliers (ADMM) \cite{boyd2011distributed} algorithm to decouple the problem into the alternated minimizations of two subproblems: a quadratic program (QP) and a spatial smoothing subproblem, where we use the QP to iteratively enforce the aggregation constraints, while solving the spatial smoothing subproblem with a recently proposed finite element technique called Spatial Spline Regression (SSR) \cite{Sanga13}. 
Moreover, our proposed algorithm not only leverages the local continuity to perform spatial sparse recovery, but is also able to take into account additional external information 
in the underlying geographical region 
to further enhance recovery performance.

We conducted extensive evaluation of the proposed algorithm in a case study of reconstructing mobile phone activity distributions of Milan, Italy from aggregated observations on base stations.
Results suggest that our algorithm achieves significant improvement, as compared to various other methods including the state-of-the-art Spatial Spline Regression (SSR) \cite{Sanga13} approach. 
Our algorithm can recover the fine-grained spatial distribution of cell phone activities in Milan only from observations on $200$ base stations, with a mean absolute percentage error of $0.309$, representing a $26.3\%$ improvement from SSR.

\section{Problem Formulation}
\label{sec:formulation}


In this section, we formally introduce the problem of recovering a spatial field from coarse aggregations observed at sparse points in the field.

Suppose the entire region of interest is modelled by an irregularly bounded domain $\Omega \subset\mathbb R^2$ that excludes the uninhabited areas such as rivers, ocean coasts, hills and so on.
Let $f(\mathbf p)$ be a real-valued function or a spatial field that models certain densities (e.g., cell phone activities) over different geographical positions $\mathbf p = (x,y)\in\Omega$.
Suppose there is a set of $m$ observation points (e.g., base stations) $B = \{B_1,\ldots, B_m\}$ sparsely distributed in $\Omega$,
where each observation point $B_i$ has a position $\mathbf p_{B_i}$ and has observed an \emph{aggregated volume} $z_i$ in the subregion $\Omega_{B_i}$ it is in charge of. 

For example, in the case of cell phone activity recovery, a mobile phone user will always connect to the closest base station (cell tower). Therefore, we have 
$
	z_i = \int_{\Omega_{B_i}}f(\mathbf p)\ud \mathbf p,
$
where the subregion $\Omega_{B_i}$ that $B_i$ represnts is given by 
\[
	\Omega_{B_i} = \{\mathbf p\in \Omega: \|\mathbf p-\mathbf p_{B_i}\|< \|\mathbf p-\mathbf p_{B_{i'}}\|, \forall B_{i'}\in B, \ i'\ne i\}.
\]
Our goal is to recover the entire spatial field $f$ of cell phone activity densities in the domain $\Omega$, based on the aggregated activities $z_1,\ldots, z_m$ observed on $m$ sparsely distributed base stations. In this case, we may call the aggregate observations $z_1,\ldots, z_m$ \emph{base station volumes}. However, such a problem is almost computationally infeasible as the continuous nature of $\Omega_{B_i}$ can hardly be handled by a personal computer.

To fix notations, suppose $\Omega$ is discretized into $n$ small grid squares $\mathbf p_1,\ldots, \mathbf p_n$, where $\mathbf p_j = (x_j,y_j)\in \Omega$, $j=1,\ldots,n$ also represents the position of square $j$'s center in $\Omega$. Without loss of generality, assume each grid square has an area of $\Delta=1$. And we have $m \ll n$, i.e, the number of aggregate observations (e.g., base station volumes) is way smaller than the number of squares to be recovered.

With the discretization of the domain, the observed volume on base station $B_i$ is given by
\begin{equation}\label{eq:BSvolumes}
	z_i = \sum_{\mathbf p_j \in \Omega_{B_i}}f(\mathbf p_j)\cdot \Delta,\quad i = 1,\ldots,m,
\end{equation}
where the subregion that $B_i$ represents is given by
\begin{equation}\label{eq:omegaBi}
	\Omega_{B_i} = \{\mathbf p_j: 1\le j\le n, \|\mathbf p_j-\mathbf p_{B_i}\|< \|\mathbf p_j-\mathbf p_{B_{i'}}\|, \ \forall i'\ne i\}.
\end{equation}
Therefore, our objective is to recover the unknown spatial field $f$, and especially the activity densities 
\[\mathbf f := (f(\mathbf p_1),\ldots, f(\mathbf p_n))^{\sf T}\] in all $n$ grid squares if the desired granularity is on a per-square level, based on the aggregated observations $z_i$ in \eqref{eq:BSvolumes}.

Note that the above problem description is not only applicable to cell phone activity density recovery, but also applies to a wide range of applications.
The nonessential difference is that each application has its own way to define each subregion $\Omega_{B_i}$, from which volume $z_i$ is aggregated.

\subsection{Constrained Spatial Smoothing Problem}

Let $\mathbf z = (z_1,\ldots,z_m)^{\sf T}$.
Since all $z_i$ are known and $\Omega_{B_i}$ can be predetermined, e.g., from \eqref{eq:omegaBi} for the cell phone activity density recovery problem, reconstructing $\mathbf f$ from \eqref{eq:BSvolumes} is essentially solving a linear equation
$
	\mathbf z = \mathbf A \mathbf f,
$
where the elements of matrix $\mathbf A\in \mathbb R^{m\times n}$ are given by 
$
A_{ij} = \left\{
	\begin{array}{ll}
		1 & \quad\text{if $\mathbf p_j \in \Omega_{B_i}$,}\\
		0 & \quad\text{otherwise.}\\
	\end{array}
\right.
$
To recover $f(\mathbf p_1),\ldots, f(\mathbf p_n)$ from $z_1,\ldots, z_m$ is apparently a sparse recovery problem, since $m\ll n$. Such a task seems infeasible, since the linear system of equations \eqref{eq:BSvolumes} is an underdetermined system which has an infinite many solutions. 

However, we may further utilize some spatial property of $f$ to make the sparse recovery problem feasible. That is, we can leverage the fact that spatial data often exhibit local correlation or local continuity within $\Omega$. For example, the cell phone activity density at a certain place critically depends on the underlying overall human population and activity at that point, e.g., downtown is more crowded than suburb residential areas, and business areas such as office buildings feature different cell phone activity patterns than leisure areas such as night clubs and restaurants. And the underlying spatial distributions of human activity density and area functionality are often slowly changing over the domain $\Omega$. 

Therefore, taking into account the local spatial continuity and non-negative property of $f$, we formulate our constrained spatial sparse recovery problem as
\begin{equation}
\label{eq:prob0}
\begin{split}
	\underset{f}{\mbox{minimize}}\quad 		& \int_\Omega(\nabla^2 f)^2\ud\mathbf p \\
	\mbox{subject to}\quad & \mathbf z = \mathbf A \mathbf f,\\
						&\mathbf f\ge 0,
\end{split}
\end{equation}
where $\nabla^2 f= \frac{\partial^2 f}{\partial x^2}+\frac{\partial^2 f}{\partial y^2}$ denotes the Laplacian of the functional $f$  to penalize the roughness of the spatial field $f$ and encourage local similarity.

Furthermore, we can also use additional external social or demographical features at each place to enhance the recovery capability.
Specifically, for each grid square $j$, let $\mathbf w_j = (w_{j1},\ldots,w_{jq})^\mathsf{T}$ be a vector of $q$ external attribute values associated with square $j$. With the additional input of external attributes, we assume that the spatial density data to be recovered in each square $j$ is given by 
\begin{equation}\label{eq:f_lr}
	f(\mathbf p_j) = f'(\mathbf p_j) + \mathbf w_j^{\sf T}\boldsymbol{\beta},
\end{equation}
where $f'(\mathbf p)$ is an underlying spatial field functional that preserves local spatial continuity, while $\mathbf w_j^{\sf T}\boldsymbol{\beta}$ is a linear regression part based on the attributes of square $\mathbf p_j$ that allows position-specific variation or jumps.

In the presence of attributes, the constrained spatial sparse recovery problem can be formulated as
\begin{equation}
\label{eq:prob1}
\begin{split}
	\underset{f', \boldsymbol{\beta}}{\mbox{minimize}}\quad 		& \int_\Omega(\nabla^2 f')^2\ud\mathbf p \\
	\mbox{subject to}
				  \quad & f(\mathbf p_j) = f'(\mathbf p_j) + \mathbf w_j^{\sf T}\boldsymbol{\beta},\quad j = 1,\ldots,n,\\
						& \mathbf z = \mathbf A \mathbf f,\\
						&\mathbf f\ge 0.
\end{split}
\end{equation}
Once the spatial field $f'$ and $\boldsymbol{\beta}$ are found, we can recover $f(\mathbf p_j)$ for all the squares using \eqref{eq:f_lr}. For example, with external attributes, the cell phone activity $f(\mathbf p_j)$ at a certain point $\mathbf p_j$ is modeled as the summation of a spatial field functional $f'(\mathbf p_j)$ and the linear regression from the attributes that permit jumps if neighboring subregions have distinct attributes and functionalities.




\section{Patched Estimation and\\ Spatial Spline Regression}
\label{sec:SSR}
In this section, we first explore some tentative solutions, and then point out their insufficiency and limitations in handling our constrained spatial sparse recovery problem.

\textbf{Patched Piece-wise Constant Estimation}.
In practice, we only know the locations of all $\mathbf p_{B_i}$'s and their corresponding aggregated volumes $z_i$'s. What is not known is the fine-grained distribution of each volume $z_i$ across $\Omega_{B_i}$, the subregion covers the observed point  $B_i$.

As a first heuristic, we can assume that the density is  distributed uniformly within each $\Omega_{B_i}$ and estimate $f(\mathbf p_j)$ as the volume $z_i$ divided by its area:
\begin{equation}\label{eq:patched}
	\bar f(\mathbf p_j) = \frac{z_i}{|\Omega_{B_i}|}, \text{ for each }\mathbf p_j \in \Omega_{B_i},
\end{equation}
where $|\Omega_{B_i}|$ is the area of $\Omega_{B_i}$.
Hence we obtain patched piece-wise constant estimation. In this paper, we may use {\it patch} to refer to $\Omega_{B_i}$, the subregion covered $B_i$. 

However, the patched estimation oversimplifies the solution, since the obtained estimates $\bar f(\mathbf p_j)$ are far from being smooth and in fact may have discontinuous jumps on the borders of patches, which will be illustrated in our evaluation. 
In practice, however, the piece-wise uniformity assumption fails: $f(\mathbf p_j)$ is not constant within a certain patch $\Omega_{B_i}$. In fact, $f(\mathbf p_j)$ should slowly change across neighboring points, as the underlying geographical and demographical characteristics also change smoothly across  regions.

\textbf{Spatial Spline Regression}.
The observations above naturally lead to the idea of using spatial smoothing techniques to smoothen the patched estimation $\bar f(\mathbf p_j)$ to remove the discontinuities and jumps. In the following, we briefly describe a recently proposed powerful smoothing technique called Spatial Spline Regression (SSR) \cite{Sanga13}. After demonstrating its usage for our particular problem, we point out the major limitations for the spatial sparse recovery problem.

Given a set of $l$ spatial data points in $\Omega$, which contains the following information: \emph{1)} the values of these $l$ points: $\{h_j\}_{j=1}^l$, \emph{2)} their positions $\{\mathbf p_j\}_{j=1}^l$, and \emph{3)} their attribute vectors $\{\mathbf w_j\}_{j=1}^l$, SSR fits a smooth spatial field $f$ by minimizing the following penalized sum of square errors \cite{Sanga13}, \cite{Ramsay02}, i.e.,
\begin{equation}\label{eq:minSSR}
	\underset{\boldsymbol{\beta}, f}{\text{minimize}}\sum_{j=1}^{l}\big(h_j- {\mathbf w}_j^{\mathsf T} \boldsymbol{\beta} - f(\mathbf p_j) \big)^2 + \lambda\int_\Omega(\nabla^2 f)^2\ud\mathbf p,
\end{equation}
where $f$ is assumed to be twice-differentiable over $\Omega$, and
$\nabla^2 f= \frac{\partial^2 f}{\partial x^2}+\frac{\partial^2 f}{\partial y^2}$ denotes the \emph{Laplacian} of $f$ to smoothen out the roughness of the spatial field $f$. The tuning parameter $\lambda$ is used to trade the smoothness of $f$ off for a better approximation to data value $h_j$.

We now briefly describe how spatial spline regression \cite{Sanga13} can solve problem \eqref{eq:minSSR} via \emph{finite element analysis} for any irregularly shaped domain $\Omega$. 
In SSR, the domain $\Omega$ is divided into small disjoint triangles, which can be done for example by the means of Delaunay triangulation \cite{hje06}. Then a polynomial function is defined on each of these triangles, such that the summation of these polynomial functions defined on different pieces closely approximates the desired spatial field $f$. It is shown in \cite{Sanga13} that the best approximation is achieved by simply solving a set of linear equations (see \cite{Sanga13} for more details).

Now we can see that if $l=n$ and we plug $h_j = \bar f(\mathbf p_j)$, $j=1,\ldots,n$ into problem \eqref{eq:minSSR}, we will get a new density surface $\hat f$ as a solution to the SSR problem \eqref{eq:minSSR} that is a smoothened approximation of the patched estimates $\bar f(\mathbf p_j)$.

However, SSR given by \eqref{eq:minSSR}
can not accommodate any constraints, and especially, does not enforce the aggregated volume constraint \eqref{eq:BSvolumes}, or equivalently, the constraint $\mathbf z = \mathbf A\mathbf f$ in  \eqref{eq:prob0}. Therefore, if we smoothen the patched estimates $\bar f(\mathbf p_j)$ out to get a smooth surface estimate $\hat f$, there is no guarantee that the estimated densities in each patch $\Omega_{B_i}$ will sum up to the observed volume $z_i$ on the point $B_i$. Violating this constraint would likely cause large density estimation errors.
\section{An ADMM Algorithm for\\ Constrained Spatial Smoothing}
\label{sec:model}

The spatial sparse recovery problem \eqref{eq:prob0} is different from \eqref{eq:minSSR} from two aspects: the loss function and the additional constraints.  As a consequence, previous SSR based on the finite element analysis can not be directly applied. More advanced algorithm needs to be developed to cope with our new loss function with constraints. 

In this section, we propose to utilize the Alternating Direction Method of Multipliers (ADMM), originated in classical paper \cite{douglas1956numerical}, to decompose our constrained optimization problem into two sub-problems that can be solved effectively by SSR and Quadratic Programming (QP) respectively. 







By introducing the following indicator function $\mathds{1}_{\mathbf{f}}$,
\begin{equation}
\mathds{1}_{\mathbf{f}} = \left\{
	\begin{array}{ll}
		0 & \quad\text{if $\mathbf{f} \geq 0$ and $\mathbf{z} = \mathbf{Af}$,}\\
		\infty & \quad\text{otherwise.}\\
	\end{array}
\right.
\end{equation}
the original problem \eqref{eq:prob0} is equivalent to
\begin{equation}
\label{eq:add-auxiliary}
\begin{split}
	\underset{f}{\mbox{minimize}} \quad 		
	& \lambda \int_\Omega(\nabla^2 f)^2\ud\mathbf p + \mathds{1}_{\mathbf{f}}
\end{split}
\end{equation}
where $\lambda$ is a parameter that controls the smoothness of $f$.

In order to split the convex optimization problem into two sub-convex problems we introduce an auxiliary variable $\mathbf g$ defined as
\begin{equation}
\mathbf g := (g(\mathbf{p}_1), \dots, g(\mathbf{p}_n))^{\sf T}
\end{equation}
The problem is transformed into the standard ADMM format,
\begin{equation}
\label{eq:add-auxiliary}
\begin{split}
	\underset{f}{\mbox{minimize}} \quad 		
	& \lambda \int_\Omega(\nabla^2 f)^2\ud\mathbf p + \mathds{1}_{\mathbf{g}}\\
	\mbox{subject to} \quad 
	& \mathbf f = \mathbf g,
\end{split}
\end{equation}
The \textit{augmented Lagrangian} for \eqref{eq:add-auxiliary} is
\begin{equation}
\label{eq:aug-Lagrangian}
\begin{split}
	{\mbox{minimize}} \quad
	\mathcal L_{\rho}(\mathbf f, \mathbf g, \boldsymbol{\alpha}) =
	& \lambda\int_\Omega(\nabla^2 f)^2\ud\mathbf p + \mathds{1}_{\mathbf{g}} \\
	& + \boldsymbol{\alpha}^{\sf T} ( \mathbf{g} -  \mathbf{f}) + \frac{\rho}{2} \| \mathbf{g} -  \mathbf{f} \|_{2}^{2},
\end{split}
\end{equation}
where
$\boldsymbol{\alpha} = (\alpha_1, ..., \alpha_n)^{\sf T}$ is the dual variable, and $\rho > 0$ is the penalty parameter in ADMM. Then the ADMM consists of the following iterations
\begin{align}
	\mathbf{f}^{k + 1} := &\ \underset{\mathbf f} {\mbox{argmin}} \ \mathcal L_{\rho} (\mathbf f, \mathbf g^{k+1}, \boldsymbol{\alpha}^k ) \label{eq:ADMM-g-1}\\
	\mathbf{g}^{k + 1} := &\ \underset{\mathbf g} {\mbox{argmin}} \ \mathcal L_{\rho} (\mathbf f^k, \mathbf g, \boldsymbol{\alpha}^k ) \label{eq:ADMM-g-2}\\
	\boldsymbol{\alpha}^{k + 1} := &\ \boldsymbol{\alpha}^{k} + \rho (\mathbf f - \mathbf g). \label{eq:ADMM-g-3}
\end{align}

For the $\mathbf f$-update step in each iteration, \eqref{eq:ADMM-g-2} is equivalent to
\begin{equation}\label{eq:admm-f}
	\underset{f}{\text{minimize}}\quad \big\| \left(\boldsymbol{\alpha}^{\sf T} + \rho \mathbf g^{\sf T}\right)/2 - \mathbf f \big\|_2^2 + \lambda\int_\Omega(\nabla^2 f)^2\ud\mathbf p,
\end{equation}
which is exactly the form of \eqref{eq:minSSR} with $h_j = \left(\alpha_j + \rho g(\mathbf p_j)\right)/2$ and $\mathbf w_j = 0$, thus can be solved efficiently by SSR. 
The penalty parameter $\lambda$ controls the smoothness of $f$ by putting little emphasis on the smoothness if it is small, and the estimated surface $f$ will be over fitted. If it is too big, the surface will be too smooth, which can cause underfitting.

For the $\mathbf g$-update step in each iteration, \eqref{eq:ADMM-g-1} is equivalent to
\begin{equation}
\label{eq:g-QR}
\begin{split}
	\underset{\mathbf g}{\mbox{minimize}} \quad 		
	& \frac{\rho}{2}\|\mathbf g\|_2^2 + (\boldsymbol{\alpha}^{\sf T} - \rho \mathbf f^{\sf T}) \mathbf g\\
	\mbox{subject to} \quad 
	& \mathbf g \geq 0,\\
	& \mathbf z= \mathbf{Ag},
\end{split}
\end{equation}
which is a convex problem that can be solve by Quadratic Programming (QP).

For the case with attributes, the algorithm does not require major changes. We just need to replace $\mathbf f$ by $\mathbf f + \mathbf W \boldsymbol{\beta}$ in \eqref{eq:admm-f}, where $\mathbf W := (\mathbf w_1, ..., \mathbf w_n)^{\sf T}$ represents the attributes and $\boldsymbol{\beta}$ is the corresponding contributions. 

Our proposed ADMM training algorithm is able to efficiently fit the spatial field and covariates for our constrained spatial sparse recovery problem. In $\mathbf g$-update step, we enforce the constraints by solving a constrained QP with no need to worry about smoothing; in $\mathbf f$-update step, we approximate the obtained $\mathbf g$ with a smooth $f$ using the SSR-based smoothing technique. In
this way, we decouple the  handling of smoothing and constraints which was not possible in pure SSR previously.

\section{Performance Evaluation}
\label{sec:simu}

\begin{figure*}[!htb]
    \centering
            \subfigure[November, $n_{\text{BS}} = 200$]{
                \includegraphics[width=1.7in]{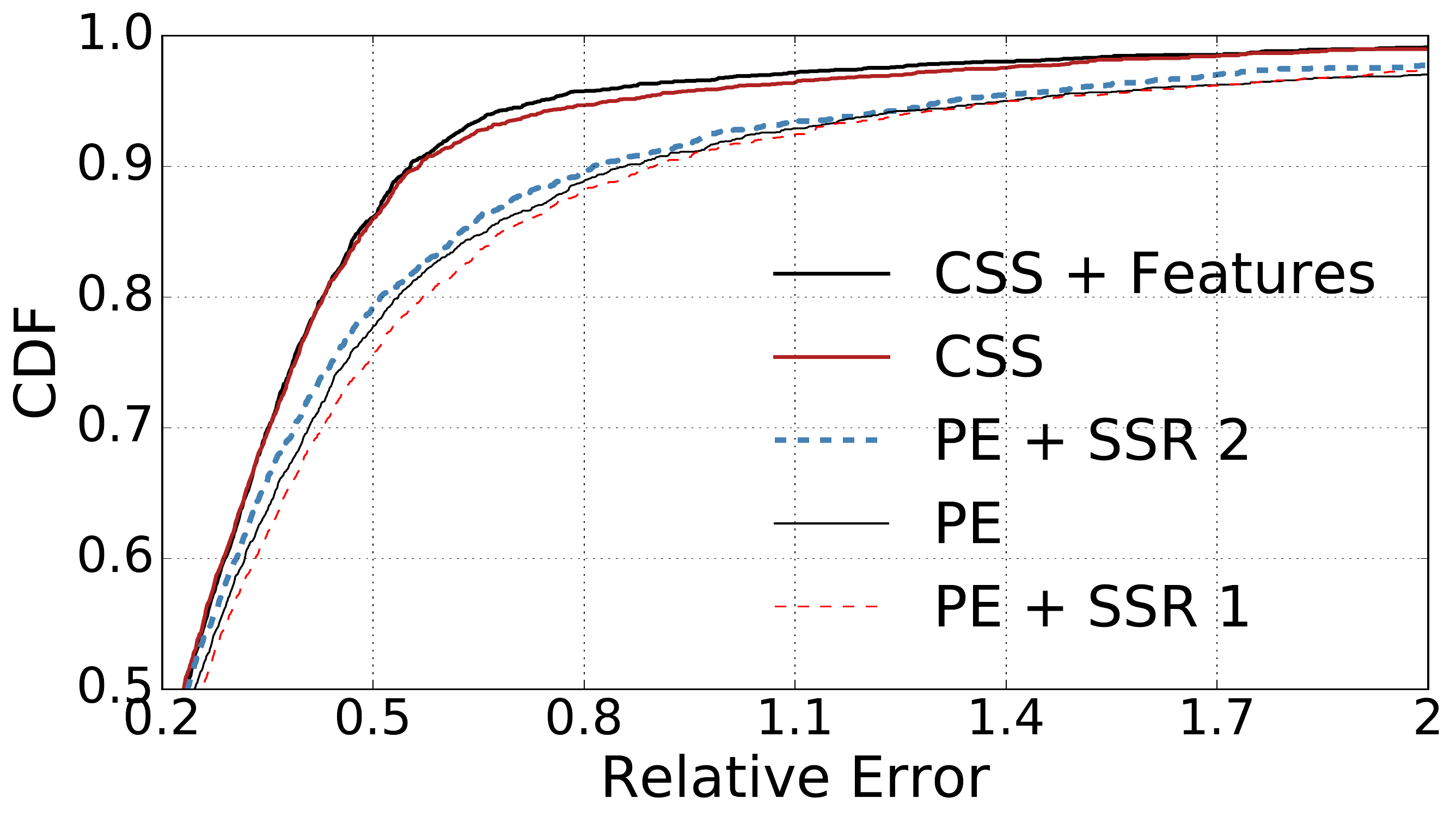}
                    \label{fig:CDF200Nov}
            }
            \hspace{-3mm}
            \subfigure[December, $n_{\text{BS}} = 200$]{
                \includegraphics[width=1.7in]{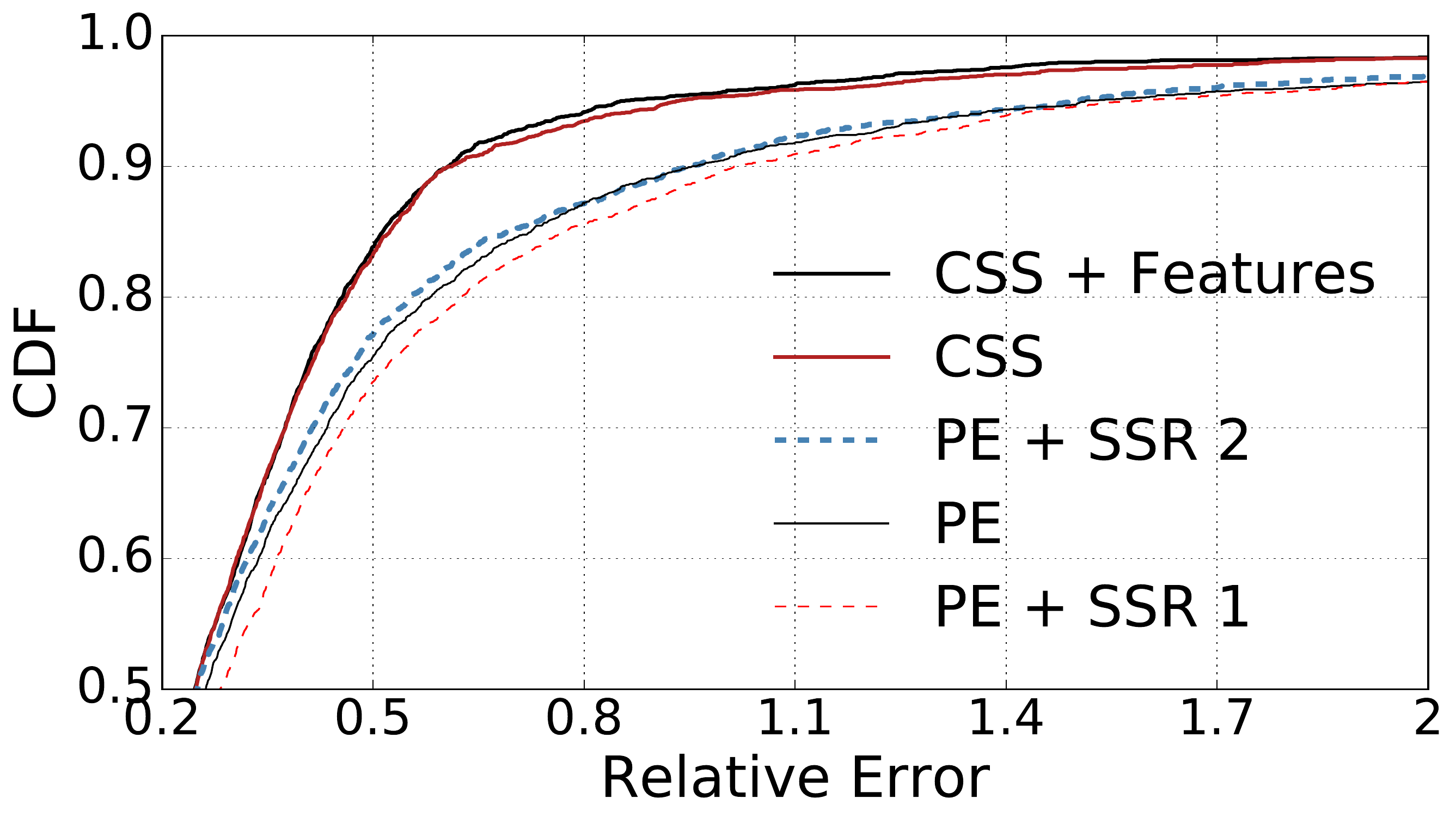}
                    \label{fig:CDF200Dec}
            }
            \hspace{-3mm}
            \subfigure[November, $n_{\text{BS}} = 100$]{
                \includegraphics[width=1.7in]{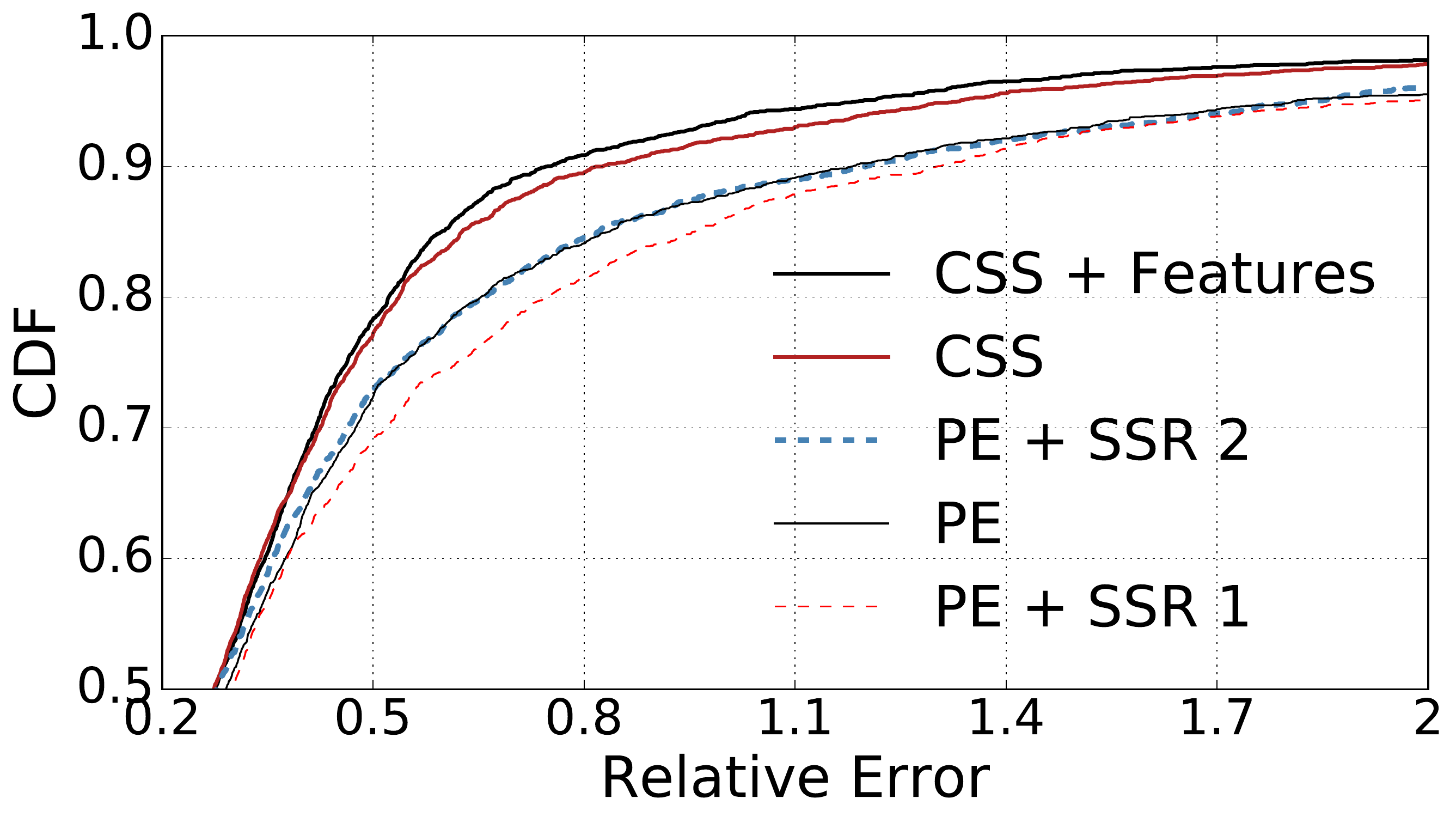}
                    \label{fig:CDF100Nov}
            }
            \hspace{-3mm}
            \subfigure[December, $n_{\text{BS}} = 100$]{
                \includegraphics[width=1.7in]{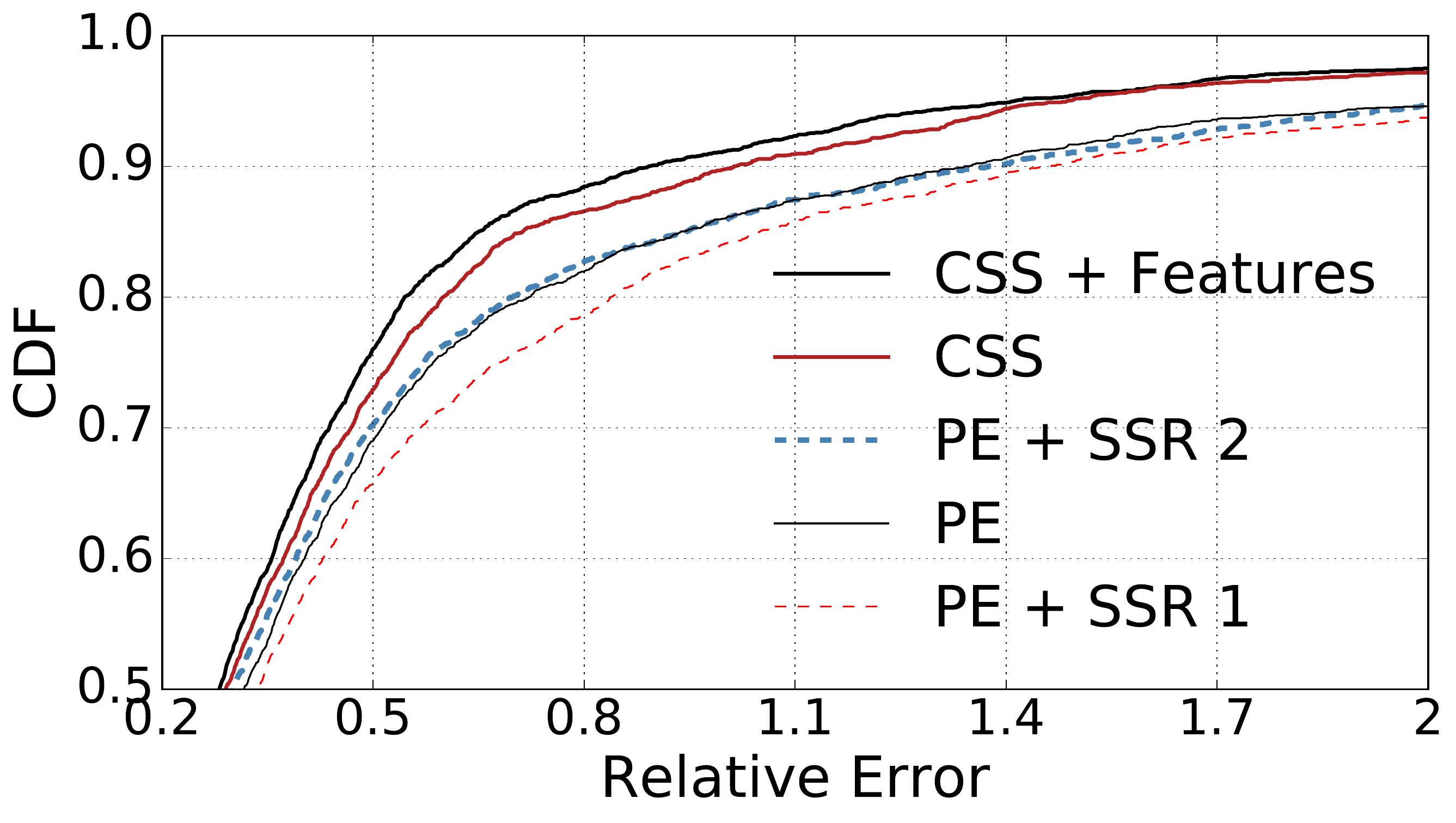}
                    \label{fig:CDF100Dec}
            }
            \vspace{-3mm}
    \caption{The comparison of the CDFs of relative errors given by different estimation methods when $n_{\text{BS}} = 200$ and $n_{\text{BS}} = 100$ for stress-testing. The legends follow the same order as the curves at relative error $= 0.5$.}
    \label{fig:compareCDF100}
\vspace{-3mm}
\end{figure*}


\begin{figure*}[t]
    \centering
            \subfigure[November, $n_{\text{BS}} = 200$]{
                \includegraphics[width=1.7in]{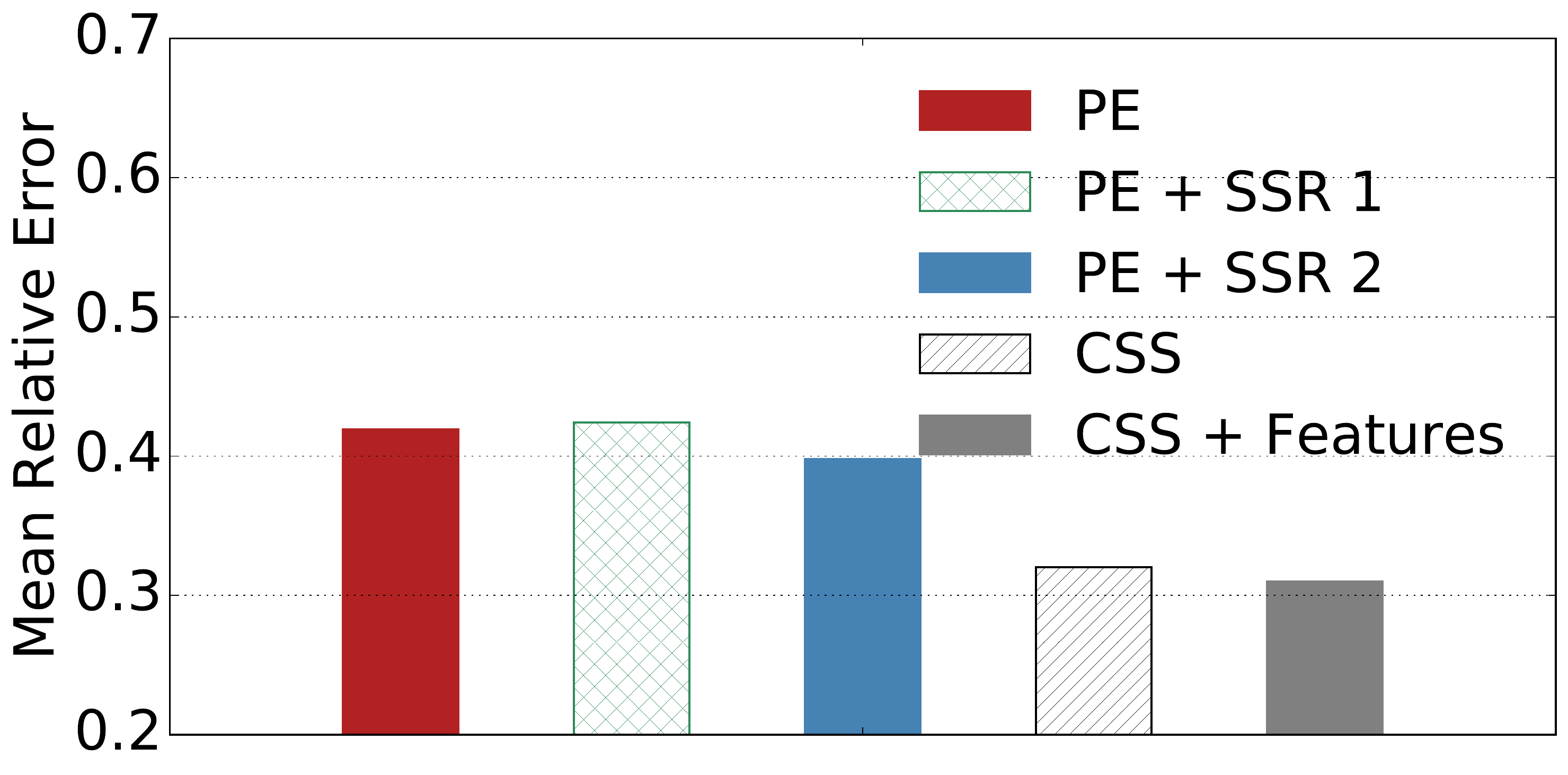}
                    \label{fig:barplot200Nov}
            }
            \hspace{-3mm}
            \subfigure[December, $n_{\text{BS}} = 200$]{
                \includegraphics[width=1.7in]{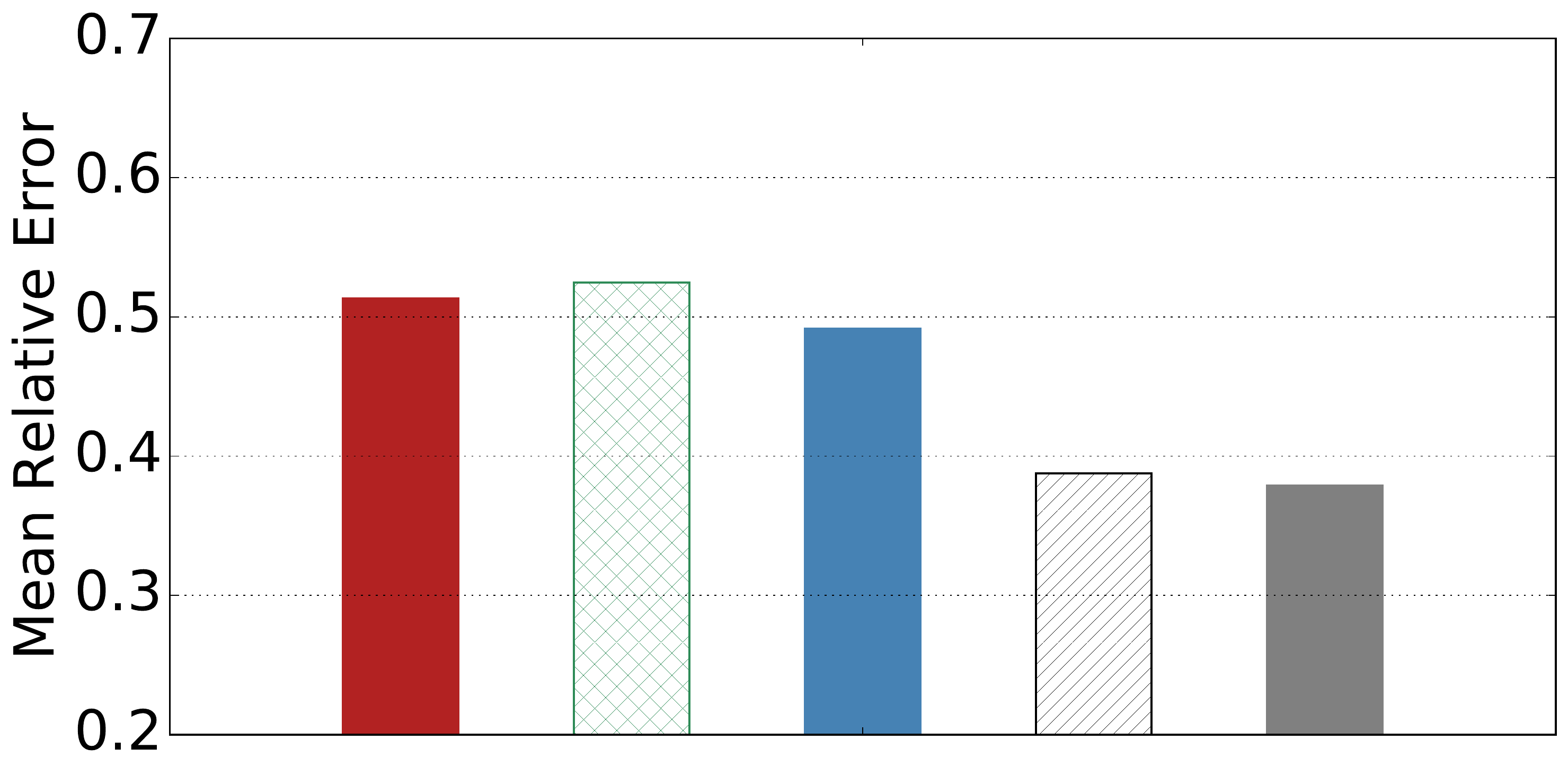}
                    \label{fig:barplot200Dec}
            }
            \hspace{-3mm}
            \subfigure[November, $n_{\text{BS}} = 100$]{
                \includegraphics[width=1.7in]{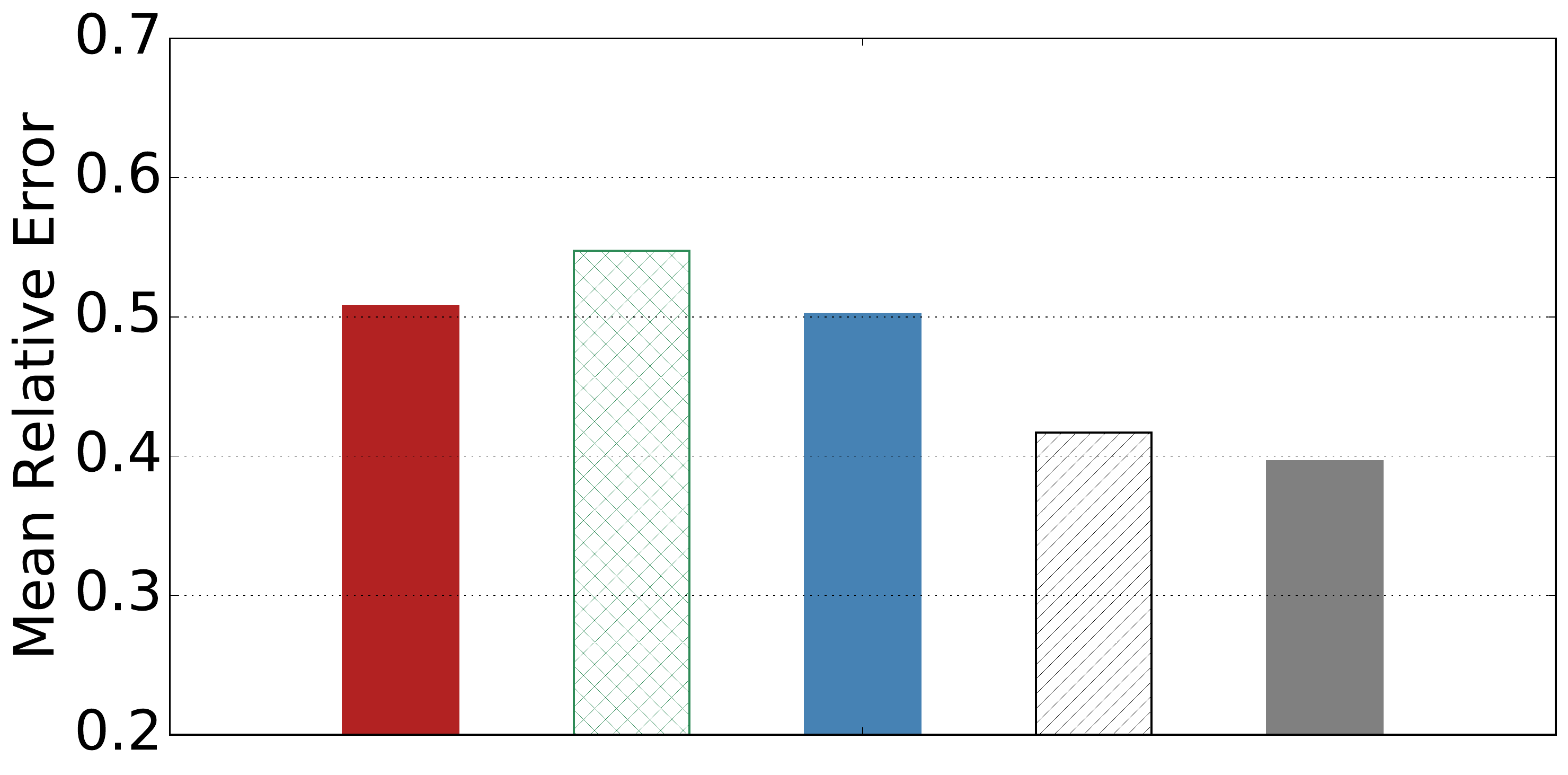}
                    \label{fig:barplot100Nov}
            }
            \hspace{-3mm}
            \subfigure[December, $n_{\text{BS}} = 100$]{
                \includegraphics[width=1.7in]{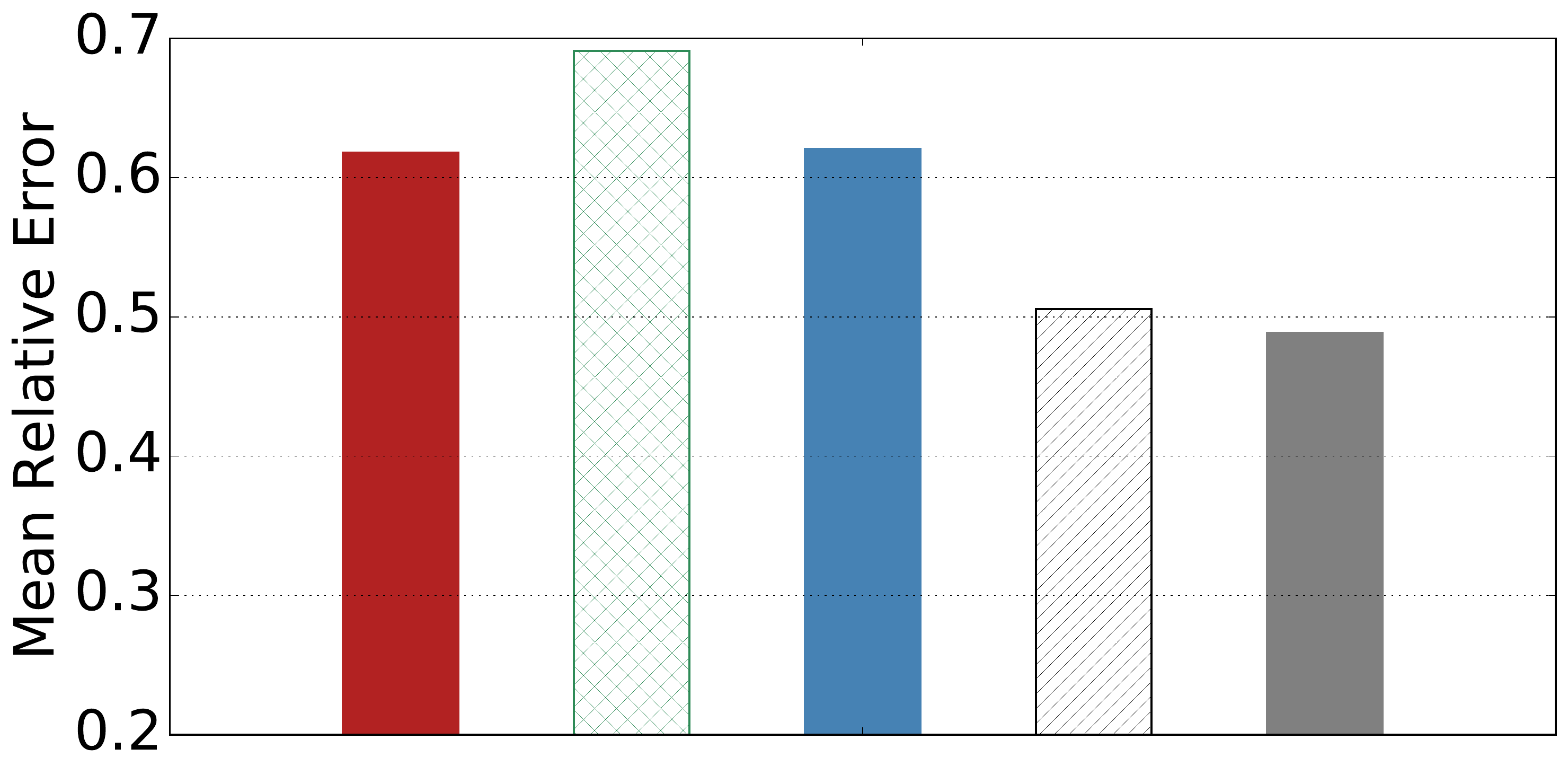}
                    \label{fig:barplot100Dec}
            }
            \vspace{-3mm}
    \caption{The comparison of the Mean Relative Error of different estimation methods when $n_{\text{BS}} = 200$ or $n_{\text{BS}} = 100$ for stress-testing. In each figure, the bars from left to right represent Patched Estimation, Patched Estimation + SSR 1, Patched Estimation + SSR 2, Constrained Spatial Smoothing, and Constrained Spatial Smoothing + Features, respectively.}
    \label{fig:compareMRE}
\vspace{-6mm}
\end{figure*}



The model in \eqref{eq:add-auxiliary} is not attached to any particular empirical problem and does not contain many implicit assumptions, it is general. However, in order to measure its performance we evaluate the model using real-world cell phone data.



The Milan Call Description Records (CDR) dataset
contains the telecommunications activity records from November 1$st$, 2013 to December 31$th$, 2013 in the city of Milan~\cite{bigdatachallenge}. In the Milan CDR dataset, the city of Milan is divided into a $100\times 100$ square grid. Each square is size of about 235m $\times$ 235m. Each activity record consists of the following entries: square ID, time-stamp of 10-minute time slot, incoming SMS activity, outgoing SMS activity, incoming call activity and outgoing call activity. The values of the four types of activities are normalized to the same scale.


Another dataset we utilized is the Milan geographical attribute dataset available from the Municipality of Milan's Open Data website \cite{barlacchi2015multi}. This dataset consists of features of central 2726 squares among the whole $10,000$ squares. The features of each square include: population, green area percentage, number of sport centers, number of universities, number of businesses, and number of bus stops.
In our empirical study, we focus on these squares to compare the performance of different algorithms.


The general \textbf{\textit{problem of recovering a spatial field from coarse aggregations observed at sparse points in
the field}} in this particular case study is reformulated into \textbf{\textit{the problem of recovering the distribution of cell phone activities over the whole 2726 square regions given that only aggregated activity observations in base stations are known}}. To study this problem, we need to further process the Milan CDR dataset.

\begin{figure}[t]
                        \centering
                        \subfigure[Distribution of BSs]{
                \includegraphics[width=1.3in]{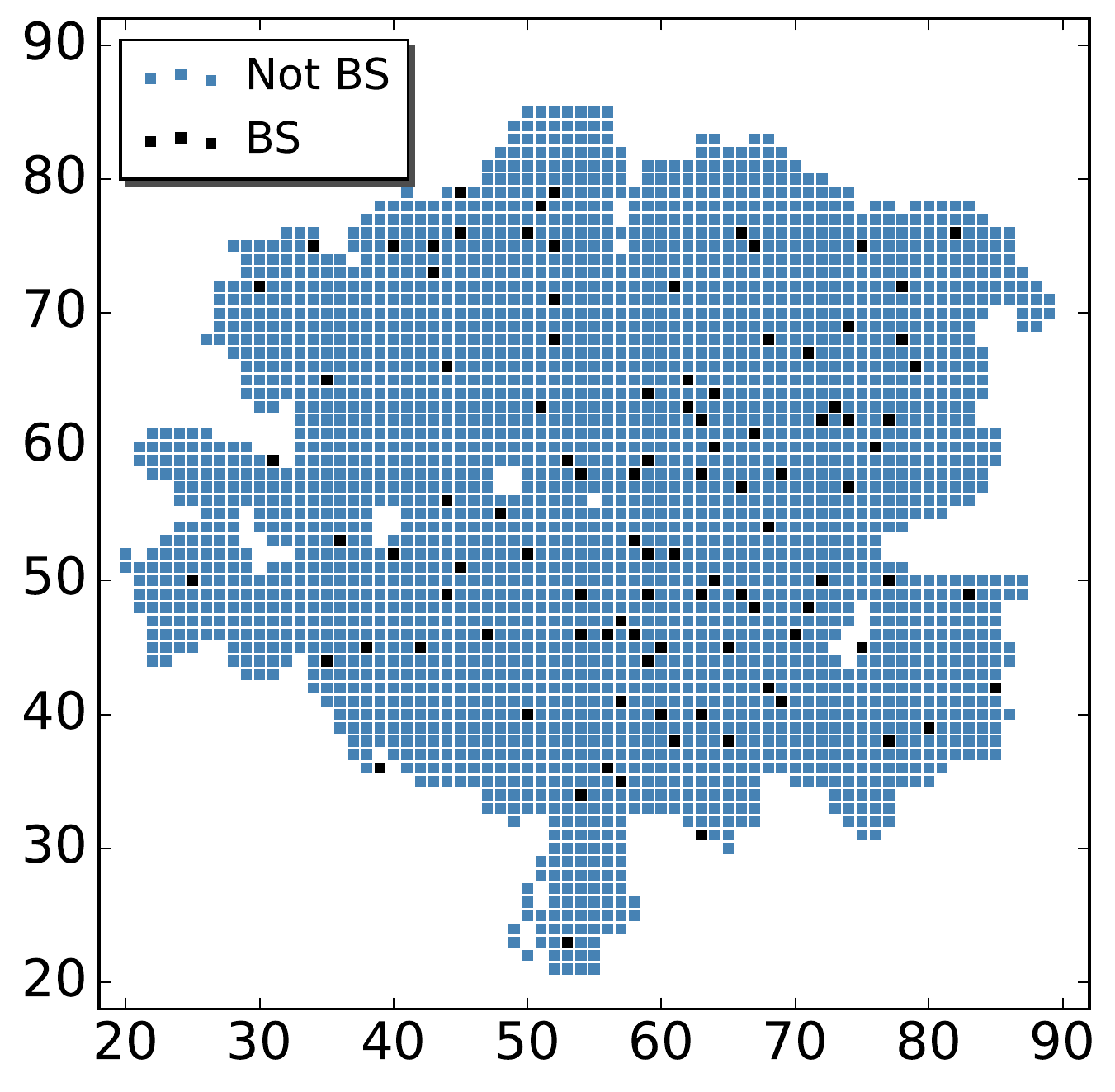}
                                \label{fig:100BS}
                        }
                \hspace{-0mm}
                        \subfigure[Areas covered by each BS]{
                \includegraphics[width=1.3in]{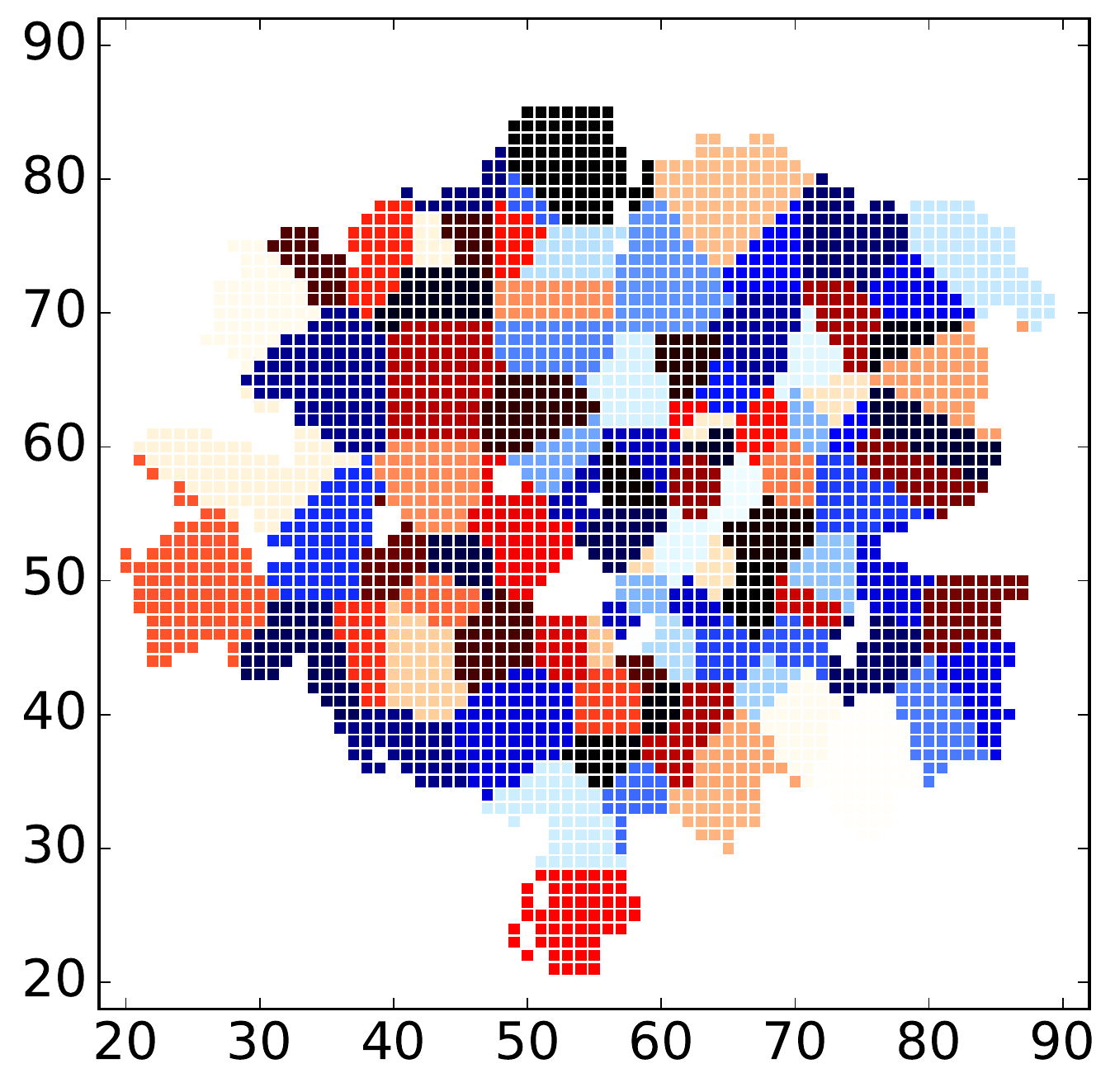}
                                \label{fig:100BSincharge}
                        }
                        \vspace{-3mm}
                \caption{(a) The geographical distribution of sampled base stations for $n_{\text{BS}} = 100$. (b) The areas that the individual base stations are responsible for, when $n_{\text{BS}} = 100$.}
                \label{fig:BSLocations}
\vspace{-6mm}
\end{figure}

\textit{First}, we sum up the four types of activities during November and December respectively to come up with the activity volume of each square during the two months. These two datasets are served as the ground-truth datasets of Milan cell phone activity distributions.
\textit{Second}, after we aggregated the two months' activities for each square, we need to set the locations of base stations (BSs). According to \cite{ratti2006mobile}, there are roughly $200$ base stations in Milan. However, the exact locations are not available. Thus, we assume the $n_{\text{BS}}$ ($n_{\text{BS}} = 200$ or $100$ for stress-test) BSs are randomly distributed according to the probability distribution
$
\Pr (\text{Set square $i$ as BS}) = f(\mathbf p_i) / \sum_{j=1}^{N}f(\mathbf p_j),
$
where $f(\mathbf p_i)$ is the cell phone activity volume in square $i$, $i=\{1, \ldots, N\}$, $N=2726$ is the number of squares we are focusing on. 
\textit{Third}, after the base station locations are sampled, the activity of each square will be assigned to its closest base station. If multiple base stations are equidistant from the square, then the activity of this square will be evenly distributed among these base stations. We then assume we only know the aggregated activities in base station squares, which is usually the true case in reality.
Fig.~\ref{fig:100BS} and Fig.~\ref{fig:100BSincharge} show the base station distributions and the region charged by each base station for $n_{\text{BS}} = 100$ respectively.


We test our proposed approach and compare it with 3 baseline methods.
\begin{itemize}
\item \textbf{Patched Estimation (PE)}:
assume cell phone activity density is distributed uniformly within each sub-region $\Omega_{B_i}$
and estimate each square's activity volume by \eqref{eq:patched}.
\item \textbf{Patched Estimation + SSR 1}: first estimate \textit{only base station} activity volumes by \eqref{eq:patched}. Use these sparse points to fit a smooth surface by running Spatial Spline Regression to obtain the estimated cell phone activity in all squares. 
\item \textbf{Patched Estimation + SSR 2}: first estimate the activity volumes of \textit{all squares} by Patched Estimation. Then use all these points to fit a smooth surface by running Spatial Spline Regression to obtain the final estimated cell phone activity in all squares.
\item \textbf{Constrained Spatial Smoothing (CSS)}: first get the initial estimation of the activity volumes of all squares by Patched Estimation, then run Constrained Spatial Smoothing algorithm to get the final activity volumes estimation of all squares.
\item \textbf{Constrained Spatial Smoothing + Features}: in this case, we incorporate the geographical features into the Constrained Spatial Smoothing algorithm.
\end{itemize}

We set the penalty parameter $\lambda = 1$ when $n_{\text{BS}} = 200$ and $\lambda = 10$ when $n_{\text{BS}} = 100$, for all methods that utilize SSR. The geographical features of Milan are only incorporated in the last algorithm described above.
We evaluate the performance by the Mean Relative Error (MRE) of the produced activity estimates for the true activity values. 

\subsection{Performance Evaluation}
\subsubsection{\bf Comparison of Different Algorithms}

We show the cumulative distribution function (CDF) of Relative Errors given by each approach in Fig.~\ref{fig:compareCDF100}. In addition, we compare the estimation's Mean Relative Errors of different approaches in Fig.~\ref{fig:compareMRE}. It is quite clear that our proposed algorithms outperform other three baseline approaches significantly in all the cases ($n_{\text{BS}} = 200$ and $n_{\text{BS}}=100$, data aggregated in November and in December). 

By comparing Patched Estimation + SSR 1 with Patched Estimation approach, we can see that using spatial smoothing based on only base station squares' observations leads to worse performance than patched estimation. This can be explained by the smoothing property of SSR and the way we set the values of base station squares. As we described, we set the activity values of base stations by averaging the total activity amount of each base station on all the squares it covers. Thus, given the activity $\frac{z_i}{|\Omega_{B_i}|}$ ($|\Omega_{B_i}|$ denotes the number of squares within region $\Omega_{B_i}$) of a base station $B_i$, the true activities of itself and its surrounding squares within region $B_i$ are distributed with a mean of $\frac{z_i}{|\Omega_{B_i}|}$. Given two base stations $B_1$ and $B_2$ that are close to each other, with aggregated activities of $z_1$ and $z_2$ respectively, the Spatial Smoothing approach will fit a smooth surface between the two base stations. Suppose $z_1 > z_2$, in this case, in overall the activities of $B_1$'s neighbour squares will be under estimated, and that of $B_2$ will be over estimated. Therefore, Patched Estimation + SSR 1's performance is not as good as Patched Estimation.

By comparing Patched Estimation + SSR 2 with Patched Estimation and Patched Estimation + SSR 1, we can observe that applying spatial smoothing on the results of patched estimation improves the performance. This proves the rationality and effectiveness of introducing smoothness into the estimated cell phone activity distribution surface.

Our proposed approaches achieves much better performance compared with the three baseline methods. By using Constrained Spatial Smoothing instead of applying Spatial Spline Regression directly, we are able to fit a smooth activity distribution while forcing it to match the observations of base station squares (the aggregated activity volumes) at the same time. By comparing Constrained Spatial Smoothing that incorporates additional features of each square with the version without features, we can see that the performance is further improved. The reason is that the heterogeneity of different locations will influence the telecommunication activity distribution, therefore making the distribution not everywhere smooth. Incorporating additional features into our model can help to explain the residuals between estimated smooth distribution and the true activity distribution, therefore further increases estimation accuracy.

The performance of different methods on December dataset is worse than on November dataset. The reason is that, there are multiple holidays during December, therefore the cell phone activities will be much more irregular than usual. 

\section{Related Work}
\label{sec:related}




The Telecom Italia Big Data Challenge dataset is a multi-source dataset that contains a variety of informations, including aggregation of telecommunication activities, news, social networks, weather, and electricity data from the city of Milan. With the important information about human activities contained in the dataset, especially the cellphone activity records, researchers utilized the data to study different problems, such as modeling human mobility patterns~\cite{gonzalez2008understanding}, population density estimation~\cite{douglass2015high}, models the spread of diseases~\cite{blondel2015survey}, modeling city ecology~\cite{bcici_mobihoc15}, etc. However, few research work has been done to estimate the spatial distribution of cellphone activity itself, despite the great value of this problem.

There are various tasks where the key problem is estimating a spatial field over a region based on observations of sampled points, such as house price estimation and population density estimation. \cite{chopra2007discovering} models the underlying surface of land desirability using kernel-based interpolation. However, it is hard to choose the form of kernel functions and tune a large number of hyper-parameters. 
Spatial Spline Regression technique is applied to the problem of population density estimation in~\cite{Sanga13}. However, in our problem, we only get the accumulated activity density in base stations, rather than real densities in each base station location. Besides, BS locations distribution is highly sparse in our case.


The fine-grained data for the distribution of the volume of calls and SMS is not usually available. A common type of data is the data collected by cell phone base stations.
Some researchers interpolate the data to obtain fine grained distributions as in \cite{ratti2006mobile}. However in \cite{ratti2006mobile} authors do not evaluate the performance of the interpolated distribution. To the best of our knowledge 
there is no extensive work done in trying to obtain optimal reconstructions of fine grained cell phone data distribution. 
We are the first to apply latest spatial functional analysis techniques to cellphone activity distribution modeling, assuming the activity densities consists of a regression part based on social or demographical statistic features and a spatial field that captures the underlying smoothness property of cellphone activities. 
\section{Concluding Remarks}
\label{sec:conclude}

In this paper, we study the problem of inferring the fine-grained spatial distribution of certain density data in a region based on the aggregate observations recorded for each of its subregions.
We propose the Constrained Spatial Smoothing (CSS) approach that exploits both the intrinsic smooth property of underlying factors and the additional features from external social or domestic statics. We further propose a training algorithm which combines the Spatial Spline Regression (SSR) technique and ADMM technique to learn our model parameters efficiently. To evaluate our algorithm and compare it with various other approaches, we run extensive evaluation based on the Milan Call Detail Records dataset provided by Telecom Italia Mobile. The simulation results on the dataset show that our algorithm significantly outperforms other baseline approaches by a great percentage. 


\bibliographystyle{IEEEtran}
\bibliography{main}

\end{document}